\documentclass{llncs}
\usepackage{amssymb}
\usepackage{amsfonts}
\usepackage{amsmath}
\usepackage{comment}
\usepackage{multirow}
\usepackage{longtable}
\usepackage{chngpage}
\usepackage{color}

\newcommand{\LL}{\mbox{$\cal L$}}
\newcommand{\PP}{\mbox{$\cal P$}}

\newcommand{\SSS}{\mbox{$\cal S$}}
\newcommand{\HH}{\mbox{$\cal H$}}

\def\R{\mathbb{R}}

\def\row#1#2{{#1}_1,\ldots ,{#1}_{#2}}
\def\brow#1#2{{\bf {#1}}_1,\ldots ,{\bf {#1}}_{#2}}

\def\R{\mathbb{R}}

\def\row#1#2{{#1}_1,\ldots ,{#1}_{#2}}
\def\brow#1#2{{\bf {#1}}_1,\ldots ,{\bf {#1}}_{#2}}

\begin{document}
\mainmatter  

\title{On the Number of Facets of Polytopes Representing  Comparative Probability Orders
}
\titlerunning{Hierarchical Simple Games: Representations and Weightedness}

\author{\bf Ilya Chevyrev, \inst{1}, 
Dominic Searles \inst{2} 
and Arkadii Slinko \inst{3}
}

\institute{Department of Mathematics, University of Auckland, Auckland, New Zealand, (\email{ichevyrev@gmail.com}) \and 
Department of Mathematics, University of Illinois at Urbana-Champaign, Urbana, USA (\email{searles2@illinois.edu})  \and Department of Mathematics, University of Auckland, Auckland,
New Zealand, (\email{a.slinko@auckland.ac.nz})
}

\maketitle

\begin{abstract}
Fine and Gill \cite{FG} introduced the geometric representation for those comparative probability orders on $n$ atoms that have an underlying probability measure. In this representation every such comparative probability order is represented by a region of a certain hyperplane arrangement. Maclagan \cite{DM} asked how many facets a polytope, which is the closure of such a region, might have. We prove that the maximal number of facets is at least $F_{n+1}$, where $F_n$ is the $n$th Fibonacci number. We conjecture that this lower bound is sharp. Our proof is combinatorial and makes use of the concept of flippable pairs introduced by Maclagan. We also obtain an upper bound which is not too far from the lower bound.
\end{abstract}
\vspace{5mm}

\bigskip


\section{Introduction}
Considering comparative probability orders from the combinatorial viewpoint, Maclagan \cite{DM}  introduced  the concept of a flippable pair of subsets.  This concept appears to be very central for the theory as it has nice algebraic and geometric characterisations. Algebraically, comparisons of subsets in flippable pairs correspond to irreducible vectors in the discrete cone associated with the comparative probability order, i.e., those vectors that cannot be split into the sum of two other vectors of the cone \cite{PF1,CCS}. Geometrically, a representable comparative probability order corresponds to a polytope in a certain arrangement of hyperplanes and flippable pairs correspond (with one exception) to those facets of the polytope which are also facets of one of the neighboring polytopes.  

Christian et al \cite{CCS} showed that in any minimal set of comparisons that define a representable comparative probability order all pairs of subsets in those comparisons are flippable.

Maclagan formulated a number of very interesting questions (see \cite[p. 295]{DM}).  In particular, she asked  how many flippable pairs a comparative probability order on $n$ atoms may have. In this paper we  show that a representable comparative probability order may have up to $F_{n+1}$ flippable pairs, which is the $(n+1)$th Fibonacci number.  We conjecture that this lower bound on the maximal number of flippable pairs is sharp. This conjecture was put forward by one of us and we call it Searles' conjecture. We provide an upper bound on the maximal number of flippable pairs in a representable comparative probability order that is not too far from $F_{n+1}$. 

Section 2 contains preliminary results and formulates Maclagan's problem. In Sections 3 and 4 we discuss Searles' conjecture in relation to Maclagan's problem and prove the aforementioned lower and upper bounds. Section 5 concludes by stating several open problems.
\section{Preliminaries}
\label{prels}

{\bf 2.1. Comparative Probability Orders.}
Given a (weak) order, that is, a reflexive, complete and transitive binary relation, $\preceq$ on a set $A$, the symbols $\prec$ and $\sim$ will, as usual, denote the corresponding (strict)  linear order and indifference relation, respectively.

\begin{definition}
Let $X$ be a finite set. A linear order $\preceq$ on $2^X$ is called a
{\em comparative probability order} on $X$ if $\emptyset\prec A$ for every
nonempty subset $A$ of $X$, and $\preceq$ satisfies de Finetti's axiom, namely
\begin{equation}
\label{deFeq}
A\preceq B \ \Longleftrightarrow \ A\cup C\,\preceq\, B\cup C ,
\end{equation}
for all $A,B,C\in 2^X$
such that $(A\cup B)\cap C=\emptyset$.
\end{definition}

As in \cite{PF1,PF2} at this stage of investigation we preclude indifferences between sets. 
For convenience, we will further suppose that $X=[n]=\{1,2,\ldots, 
n\}$ and denote the set of all comparative probability orders on 
$2^{[n]}$ by $\PP_n$.  

If we have a probability measure ${\bf p}=(\row pn)$ on $X$,
where $p_i$ is the probability of $i$, then we know the
probability $p(A)$ of every event $A$ which is given by $p(A)=\sum_{i\in A}p_i$.
We may now define an order $\preceq_{\bf p}$ on $2^X$ by
\[
A\preceq_{\bf p} B \quad \mbox{if and only if}\quad  p(A)\le p(B).
\]
If the probabilities of all events are
different, then $\preceq_{\bf p}$ is a comparative probability order on
$X$. Any such order is 
called {\em (additively) representable}.  The set of representable 
orders is denoted by $\LL_n$. It is known \cite{KPS} that $\LL_n$ is 
strictly contained in $\PP_n$ for all $n\ge 5$.

Since a representable comparative probability order does not have a unique probability measure representing it but a class of them, any representable comparative probability order can be viewed as a credal set (a closed and convex set of probability measures, see, e.g.,  \cite{L1}) of a very special type. We will return to this interpretation slightly later.\par

As in \cite{PF1,PF2}, it is often 
convenient to assume that
$
1\prec 2\prec \ldots \prec n.
$
This reduces the number of possible orders under consideration by a factor of $n!$. The set of all comparative probability orders on $[n]$ that satisfy this condition will be denoted by $\PP_n^*$, and the set of all such representable comparative probability orders on $[n]$ will be denoted by~$\LL_n^*$.\par

We can also define a representable comparative probability order by any vector of positive utilities ${\bf u}=(\row un)$ by
\[
A\preceq_{\bf u} B \quad \mbox{if and only if}\quad  \sum_{i\in A}u_i\le \sum_{i\in B}u_i.
\]
We do not get anything new since this will be the order $\preceq_{\bf p}$ for the measure ${\bf p}=\frac{1}{S}{\bf u}$, where $S=\sum_{i=1}^nu_i$. However, sometimes it is convenient to have the coordinates of ${\bf u}$ integers. In this case we will call $u(A)=\sum_{i\in A}u_i$ the {\em utility} of $A$.

Kraft et al \cite{KPS} gave necessary and sufficient conditions for a comparative probability order to be representable. They are not so easy to formulate and they have appeared in the literature in various forms (see, e.g., \cite{DS,PF1}). The easiest way to formulate them is through the concept of a {\em trading transform} introduced in \cite{TZ}.

\begin{definition}
A sequence of subsets $(\row Ak; \row Bk)$ of $[n]$ of even length $2k$ is said to be a trading transform of length $k$ if for every $i\in [n]$
$$
\left|\{j\mid i\in A_j\}\right|=\left|\{j\mid i\in B_j\}\right|.
$$
In other words, sets $\row Ak$ can be converted into $\row Bk$ by rearranging their elements. 
\end{definition}

Now the result of \cite{KPS} can be reformulated as follows. 

\begin{theorem}[Kraft-Pratt-Seidenberg]
\label{KPStheorem}
A comparative probability order $\preceq $ is representable if and only if for no $k$ there exist pairs $A_i\prec B_i$, $i=1,2,\ldots, k$ such that 
$
(\row Ak; \row Bk)
$
is a trading transform of length~$k$.
\end{theorem}

\par \medskip

{\bf 2.2. Discrete Cones.}
To every linear order $\preceq$~$\in \PP_n^*$, there corresponds a {\em discrete cone} $C(\preceq)$ in $T^n$, where $T=\{-1,0,1\}$ (as defined in \cite{AK,PF1}).

\begin{definition}
A subset ${\cal C}\subseteq T^n$ is said to be a discrete cone if the 
following properties hold:
\begin{enumerate}
\item[{\rm D1.}] $\{ {\bf e}_1, {\bf e}_2,\ldots, {\bf 
e}_{n}
\}\subseteq {\cal C}$,
where $\{{\bf e}_1,\ldots,{\bf e}_n\}$ is the standard basis of $\R^n$,
\item[{\rm D2.}]  for every ${\bf x}\in T^n$, exactly one vector of the set $\{-{\bf x},{\bf x}\}$ belongs to $\cal C$,
\item[{\rm D3.}] ${\bf x}+{\bf y}\in C$ whenever ${\bf x},{\bf y}\in {\cal C}$ and
${\bf x}+{\bf y}\in T^n$.
\end{enumerate}
\end{definition}
We note that in \cite{PF1} Fishburn requires ${\bf 0}\notin {\cal C}$  because
his orders are anti-reflexive. In our case, condition D2 implies 
${\bf 0}\in {\cal C}$.
\medskip

For each subset $A\subseteq X$ we define the characteristic vector 
$\chi_A$ of this subset by setting $\chi_A(i)=1$, if $ i\in A$, and $\chi_A(i)=0$, if $ i\notin A$.
Given a comparative probability order $\preceq $ on $X$,
we define the characteristic vector
$\chi(A,B)=\chi_B -\chi_A\in
T^n$  for every possible pair $(A,B) $ such that $A\preceq B$.
The set of all characteristic
vectors $\chi(A,B)$  is denoted 
by $C(\preceq)$. The two axioms of comparative probability guarantee 
that $C(\preceq)$ is a discrete cone (see \cite[Lemma~2.1]{PF1}).

\medskip

{\bf 2.3. Critical and flippable pairs.} Not all relations $A\prec B$ for pairs of subsets $(A,B)$ in a comparative probability order are equally informative. Some of these may be implied by others through transitivity or de Finetti's axiom. This is certainly true for any pair consisting of two nonadjacent sets or two sets with nonempty intersection.

\begin{definition}
Let $A$ and $B$ be disjoint subsets of $[n]$. The pair $(A,B)$ is said to be {\em critical} for $\preceq $ if $A\prec B$ and $(A,B)$ are adjacent, i.e., there is no $C\subseteq [n]$ for which $A\prec C\prec B$.
\end{definition}

In the above definition we follow Fishburn \cite{PF3}, while Maclagan \cite{DM} calls such pairs {\em primitive}.

It is known \cite{DM} that, if $\preceq$ and $\preceq'$ are distinct comparative probability orders, then there exists a critical pair $(A,B)$ for $\preceq$ such that $B \prec' A$.
This shows that critical pairs are of interest due to the fact that they define orders. But there are even more interesting pairs.

\begin{definition}
A critical pair $(A,B)$ is said to be {\em flippable\/} for $\preceq $ if  for every $D\subseteq [n]$, disjoint from $A\cup B$, the pair $(A\cup D,B\cup D)$ is adjacent in $\preceq$.
\end{definition}

We note that the set of flippable pairs is not empty, since the central pair of any comparative probability order is flippable \cite{KPS}. Indeed, this consists of a certain set $A$ and its complement  $A^c=X\setminus A$, and there is no $D$ which has empty intersection with both of these sets. It is not known whether this can be the {\em only\/} flippable pair of the order.\par
\medskip

Suppose now that a pair $(A,B)$ is flippable for a comparative probability order $\preceq  $, and $A\ne \emptyset$. Then reversing each comparison $A\cup D \prec B\cup D$ (to $B\cup D \prec A\cup D$), we will obtain a new comparative probability order $\preceq  '$, since the de Finetti axiom (\ref{deFeq}) will still be satisfied.  We say that $\preceq ' $ is obtained from $\preceq  $ by {\em flipping over\/} $A\prec B$. The orders $\preceq $ and $\preceq '$ are called {\em flip-related}. This flip relation turns $\PP_n$ into a graph which we will denote ${\cal G}_n$. 


\begin{definition}
An element ${\bf w}$ of the cone ${\cal C}$ is said to be {\em reducible\/} if there exist two other vectors
${\bf u}, {\bf v}\in {\cal C}$ such that ${\bf w}={\bf u} + {\bf v}$, and {\em irreducible\/}  otherwise. The set of all irreducible elements of ${\cal C}$ will be denoted as $\text{Irr}({\cal C})$.
\end{definition}

\begin{theorem}[\cite{DM,CCS}]
\label{flippable-irreducible}
A pair $(A,B)$ of disjoint subsets is flippable for $\preceq $ if and only if the corresponding characteristic vector  $\chi (A,B)$ is irreducible in ${\cal C}(\preceq )$.
So the cardinality $|\text{Irr}({\cal C}(\preceq))|$ is the total number of flippable pairs in $\preceq$.
\end{theorem}

Flippable pairs uniquely define a representable order but this does not hold for nonrepresentable orders \cite{CCS}.
\medskip

As we know, the flip relation turns $\PP_n$ into a graph ${\cal G}_n$. Let $\preceq $ and $\preceq '$ be two comparative probability orders which are connected by an edge in this graph (and so are flip-related). We say that $\preceq $ and $\preceq '$ are  {\em friendly\/}  if they are either both representable or both nonrepresentable. \par\medskip

{\bf 2.4. Geometric Representation of Representable  Orders and Maclagan's Problem.}
Let $A,B\subseteq [n]$ be disjoint subsets, of which at least one 
is nonempty. Let $H(A,B)$ be a hyperplane  consisting of all points ${\bf x}\in \R^n$ satisfying the equation
\begin{equation*}
\sum_{a\in A}x_a-\sum_{b\in B}x_b=0.
\end{equation*}
We denote the corresponding 
hyperplane arrangement by ${\cal A}_n$.
Also let $J$ be the hyperplane $x_1+x_2+\ldots 
+x_n=1,$ and let ${\cal H}_n={\cal A}_n^J$ be the induced hyperplane arrangement. Fine and Gill \cite{FG} showed that the regions of ${\cal H}_n$ in the 
positive orthant $\R^n_+$ of $\R^n$ correspond to representable orders from ${\PP}_n$. 

Now we can see what is special in the credal sets that correspond to comparative probability orders. They are not only convex, as credal sets must be, but they are in fact interiors of polytopes.
When in the future we refer to a region of this hyperplane arrangement we will refer to the polytope which is the closure of that region. This will invite no confusion.

\begin{problem}[Maclagan \cite{DM}]
\label{prob1}
What is an upper bound for the number of representable neighbors for a representable comparative probability order on $n$ atoms? In other words, how many facets can regions of ${\cal H}_n$ have?
\end{problem}
The maximal number of facets of regions of ${\cal H}_n$ we will call the $n$th {\em Maclagan number} and denote $M(n)$, while the maximal number of flippable pairs for a representable order on $n$ atoms will be denoted $m(n)$. In this paper we provide bounds on these functions, some of which we suspect to be sharp. 
  It is clearly sufficient to solve Maclagan's problem (Problem~\ref{prob1}) for comparative probability orders in $\LL_n^{*}$. \par\medskip

The main combinatorial  tool for calculating or estimating $M(n)$ is the following semi-obvious proposition. 

\begin{proposition}[\cite{DM,CCS}]
\label{mainprop}
Let $\preceq $ be a representable comparative probability order, and let $P$ be the corresponding convex polytope, which is a region of the hyperplane arrangement ${\cal H}_n$. Then the number of  facets of $P$ equals the number of representable comparative probability  orders that are flip-related to $\preceq$ (plus one if the pair $\emptyset \prec 1$ is flippable).
\end{proposition} 

\begin{corollary}
\label{Mandm}
$M(n)\le m(n)$.
\end{corollary}

\begin{proof}
From the proposition it follows that $M(n)$ cannot be greater that $m(n)$. However theoretically it can be smaller since not all flips of the representable comparative probability order that has the maximal number of flips may be friendly.
\end{proof}

It is worth noting that the minimal number of facets of a region in ${\cal H}_n$ is known and equal to $n$ \cite{CS,CCS}.

\section{The Lower Bound}

It is known that $M(3)=m(3)=3$ and $M(4)=m(4)=5$ \cite{CCS}.
Computations in {\sc Magma\/} {\cite{CCS}} show that $5\le |{\rm Irr}({\cal C})|\le 8$ for $n=5$ and $5\le |{\rm Irr}({\cal C})|\le 13$ for $n=6$
with all intermediate values being attainable for both values of $n$.
It was also observed that for $n=5$ and $n=6$, all comparative probability orders with the largest possible number of flips (namely 8 for $n=5$, and 13 for $n=6$) are representable, and all of their flips are friendly. This means that $M(5)=m(5)=8$ and $M(6)=m(6)=13$.

Searles noticed that the four known values are  Fibonacci numbers, i.e., belong to the sequence defined by $F_1=F_2=1$ and $F_{n+2}=F_{n+1}+F_n$. He conjectured that \par
\medskip

\noindent {\bf Conjecture (Searles, 2007)}
The maximal number of facets of regions of $\HH_n$ is equal to the maximal cardinality of $\text{Irr}({\cal C(\preceq)})$ for  $\preceq \in \LL^{*}_n$, and equal to the Fibonacci number $F_{n+1}$ or, alternatively, $M(n)=m(n)=F_{n+1}$.\par\medskip

The first part of this conjecture will be proved if we show that for some representative comparative probability order $\preceq $, for which  $|\text{Irr}({\cal C(\preceq)})|$ is maximal, all flips of $\preceq $ are friendly. The existence of such an order was checked for all $n\le 12$. \par\medskip

In this section we prove that $M(n)\ge F_{n+1}$. To this end we prove

\begin{theorem}
\label{Dtheorem}
In ${\cal P}_n$ there exists a representable comparative probability order which (a) has $ F_{n+1}$ flippable pairs and (b)  whose flips are all friendly.
\end{theorem}

The proof will be split into several observations. Let us introduce the following notation first. Let ${\bf u}=(\row un)$ be a vector such that $0<u_1<\ldots<u_n$ and $q>0$ be a number such that $u_j<q<u_{j+1}$ for some $j=0,1,2,\ldots,n$ (we assume that $u_0 = 0$ and $u_{n+1}=\infty$). In this case we set $({\bf u},q)$ to be the vector of $\R^{n+1}$ such that
\[
({\bf u},q)=(u_1,\ldots,u_j,q,u_{j+1},\ldots, u_n).
\]
We also denote ${\bf \ell}_n=(1,2,4,\ldots,2^{n-1})$ and $2{\bf \ell}_n=(2,4,8,\ldots,2^{n})$. We start with an easy and well-known observation.

\begin{proposition}
\label{pr2}
$\preceq_{{\bf\ell}_n}$ is the lexicographic order on $2^{[n]}$.  The utilities of subsets from $2^{[n]}$ cover the whole range of integers between 0 and $2^n-1$ and the utilities of any two consecutive  subsets in it differ by~$1$. 
\end{proposition}

\begin{proof} This is equivalent to every natural number possessing a unique binary representation. We leave the verification to the reader.\end{proof}

\begin{proposition}
\label{pr3}
Let $q$ be an odd positive integer smaller than $2^{n}$ and ${\bf m}=(2{\bf\ell}_n,q)$. Consider the order $\preceq_{\bf m}$ on $2^{[n+1]}$. Then the difference between the utilities of any two consecutive subsets in this order is not greater than~$2$.
\end{proposition}

\begin{proof}
Suppose $2^{j-1} \leq q<2^{j}$, that is, $q$ is the utility of $j$ in $\preceq_{\bf m}$.  By Proposition~\ref{pr2} the utilities of the subsets from $[n+1]\setminus \{j\}$ cover the range of even values from 0 to $2^{n+1}-2$.  Suppose $B$ is a subset in $\preceq_{\bf m}$, where $B \neq \emptyset$. If $u(B) \leq 2^{n+1}-2$, then by Proposition 2 there exists a subset $A$ such that $0 < u(B) - u(A) \leq 2$. If $u(B)> 2^{n+1}-2$, then we must have $j \in B$, and since $u(j) < 2^n$, $B'=B \setminus \{j\} \neq \emptyset$. As $j \notin B'$ we have $u(B') \leq 2^{n+1}-2$, and so by Proposition 2 there exists $A' \subseteq [n+1] \setminus \{j\}$ such that $A' \prec_m  B'$ and $u(B')-u(A') =2$. Then adding $j$ to both subsets we obtain $u(B) - u(A) = 2$ for $A = A' \cup \{j\}$. Therefore, for any nonempty $B$ in $\preceq_{\bf m}$, there exists a subset $A$ such that $0 < u(B) - u(A) \leq 2$, and so for any adjacent pair $(C,D)$ of subsets, we have $0 < u(D) - u(C) \leq 2$.
\end{proof}

Let us denote by $\SSS_{n+1}$ the class of orders on $X=\{1,2,\ldots,n+1\}$ of type $\preceq_{\bf m}$, where ${\bf m}=(2{\bf\ell}_n,q)$ for some odd $0 < q< 2^{n}$. And let $j$ denote the number such that $2^{j-1} \leq q<2^{j}$. Obviously, $j<n+1$. By $\row u{n+1}$ we will denote the respective utilities of elements of $X$, that is ${\bf m}=(\row u{n+1})$.

\begin{proposition}
\label{pr4}
From the position at which the subset $\{j\}$ appears in the order $\preceq_{\bf m}$ until the position after which all subsets contain $j$, subsets not containing $j$ alternate with those containing $j$, with the difference in utilities for any two consecutive terms being $1$. 
\end{proposition}

\begin{proof}
All subsets not containing $j$ have even utility and all those containing $j$ have odd utility. If we consider these two sequences separately, by Proposition~\ref{pr2} the difference of utilities of neighboring terms in each sequence will be equal to 2. Hence they have to alternate in $\preceq_{\bf m}$.
\end{proof}

\begin{lemma}
\label{lemma_abc}
Let  $\preceq_{\bf m}$ be an order from the class $\SSS_{n+1}$ and let $(A,B)$ be a critical pair for $\preceq_{\bf m}$. Then the following conditions are equivalent:
\begin{enumerate}
\item[\rm (a)] $(A,B)$ is flippable;
\item[\rm (b)] either $A$ or $B$ contains $j$ but not both;
\item[\rm (c)] $u(B)-u(A)=1$.
\end{enumerate}
\end{lemma}

\begin{proof}
(a) $\Longrightarrow $ (b): Suppose $(A,B)$ is flippable. As $(A,B)$ is critical, it is impossible for $A$ and $B$ each to contain $j$ as $A\cap B=\emptyset$. We only have to prove that it is impossible for both of them not to contain $j$. If $j\notin A$ and $j\notin B$, then $u(A)+2=u(B)$. Since the pair is critical, by Proposition~\ref{pr4} both $A$ and $B$ appear in the order earlier than $\{j\}$. Hence $u(A)+2=u(B) < u(j)$. Then, in particular, $u(A)<u(B)<u(n+1)=2^n$, hence neither $A$ nor $B$ contains $n+1$. But then for $A'=A\cup \{n+1\}$ and $B'=B\cup \{n+1\}$ we have $u(j)<u(A')<u(B')$. Both $A'$ and $B'$ do not contain $j$, hence they are in the alternating part of the order, and since $u(B')-u(A')=2$, they cannot be consecutive terms. As $(A,B)$ is flippable, this is impossible, which proves that either $A$ or $B$ contains~$j$. 

(b) $\Longrightarrow $ (c): This follows from Proposition~\ref{pr4}.

(c) $\Longrightarrow $ (a):  This is true not only for orders from our class, but also for all orders defined by integer utility vectors. Indeed, if $u(B)-u(A)=1$, then for any $C\cap (A\cup B)=\emptyset$ we have  $u(B\cup C)-u(A\cup C)=1$,  and so $A\cup C$ and $B\cup C$ are consecutive.
\end{proof}

Up to now, the values of $q$ and $j$ did not matter. Now we will try to maximise the number of flippable pairs in $\preceq_{\bf m}$, so we will need to choose them carefully. It should come as no surprise that the optimal choice of $j$ and $q$ will depend on $n$, so we will talk about $j_n$ and $q_n$ now. For the rest of the proof we will set 
\begin{equation}
\label{numberq}
j_n=n-1,\qquad q_n=\frac{(-1)^{n+1}+2^n}{3}.
\end{equation}
An equivalent way of defining $q_n$ would be by the recurrence relation
\begin{equation}
\label{recrel}
q_n=q_{n-1}+2q_{n-2}
\end{equation} 
with the initial values $q_3=3$, $q_4=5$. We also note:
\begin{proposition}
\label{pr5}
$
q_n\equiv 2+(-1)^{n+1} \pmod4.
$
\end{proposition}

\begin{proof}
Easy induction using (\ref{recrel}).
\end{proof}

Let us now consider a flippable pair $(A,B)$ for $\preceq_{\bf m}$, where ${\bf m}=(2{\bf\ell}_n,q_n)$. Since $j_n=n-1$, we have either $A=A'\cup\{n-1\}$ or $B=B'\cup\{n-1\}$ but not both. In the first case, $(A',B)$ is a pair of nonintersecting subsets from the lexicographic order induced by $2{\bf\ell}_n$ on $ [n+1]\setminus \{n-1\}$ with $u(B)-u(A')=q_n+1$. In the second, $(B',A)$ is a pair of nonintersecting subsets from the same lexicographic order with $u(A)-u(B')=q_n-1$.

As $ [n+1]\setminus \{n-1\}$ can be identified with $[n]$, we let $g_n$ be the number of pairs $(A,B)$ in the lexicographic order $\preceq_{2{\bf \ell}_n}$ on $n$ atoms with $u(B)-u(A)=q_n+1$, and let $h_n$ be the number of pairs $(A,B)$ in the same order with $u(B)-u(A)=q_n-1$. What we have proved is the following:

\begin{lemma}
\label{g_n+h_n}
Let ${\bf m}=(2{\bf\ell}_n,q_n)$. Then the number of flippable pairs in  $\preceq_{\bf m}$ is $g_n+h_n$.
\end{lemma}

This reduces our calculations to a rather understandable lexicographic order $\preceq_{2{\bf \ell}_n}$.\par\medskip

For convenience we will  denote $q_n^+=q_n+1$ and $q_n^-=q_n-1$. We note that Proposition~\ref{pr5} implies
\begin{proposition}
\label{pr6}
$
q_n^-\equiv 1+(-1)^{n+1} \pmod4$, and  $q_n^+\equiv 3+(-1)^{n+1} \pmod4.
$
In particular, if $n$ is even, $q^-_n\equiv 0 \pmod4$ and $q^+_n\equiv 2 \pmod4$ and if $n$ is odd, $q^-_n\equiv 2 \pmod4$ and $q^+_n\equiv 0 \pmod4$.
\end{proposition}

A direct calculation also shows that the following equations hold:
 
 \begin{proposition}
\label{pr7}
\begin{eqnarray}
\label{-2-}
q_{n+1}^-&=&2q^-_n \qquad\ \  \text{for all odd $n\ge 3$},\\
\label{-3-}
q_{n+1}^-&=&2q^-_n+2 \quad \text{for all even $n\ge 4$},\\
\label{-4-}
q_{n+1}^+&=&2q^+_n-2 \quad \text{for all odd $n\ge 3$},\\
\label{-5-}
q_{n+1}^+&=&2q^+_n \qquad\ \  \text{for all even $n\ge 4$}.
\end{eqnarray}
\end{proposition}

\begin{lemma}
\label{ghrecrel}
 The following recurrence relations hold:  for any odd  $n\ge 3$
\begin{eqnarray*}
g_{n+1}=g_n+h_n,\qquad
h_{n+1}=h_n,
\end{eqnarray*}
and for any even $n\ge 4$
\begin{eqnarray*}
g_{n+1}=g_n,\qquad 
h_{n+1}=g_n+h_n.
\end{eqnarray*}
\end{lemma}

\begin{proof}
Firstly we assume that $n$ is odd. Then $n+1$ is even. We know from (\ref{-2-}) that $q_{n+1}^-=2q^-_n$. Given any nonintersecting pair $(A, B)$ of subsets in $[n]$, with $A \prec_{2{\bf \ell}_n} B$ and $u(B)-u(A)=q^-_n$, we may shift both subsets to the right, replacing each element $i$ in them with the element $i+1$, to obtain a nonintersecting pair $(\overline{A},\overline{B})$ of subsets in $[n+1]$, where $\overline{A}$ precedes $\overline{B}$ in $\preceq_{2{\bf \ell}_{n+1}}$.  This procedure of shifting doubles the difference in utilities, so $u(\overline{B})-u(\overline{A})=2q^{-}_n=q_{n+1}^{-}$. This proves $h_{n+1}\ge h_n$. Moreover, by (\ref{-2-}) and Proposition~\ref{pr6}, $q^-_{n+1}\equiv 0 \pmod4$, hence no nonintersecting pair $(C, D)$ in $\preceq_{2{\bf \ell}_{n+1}}$ with difference of utilities $q^-_{n+1}$ can include~$1$, either in $C$ or in $D$, as $u_1 = 2$. Therefore $C=\overline{A}$ and $D=\overline{B}$ for some nonintersecting pair $(A, B)$ in $[n]$ with $u(B)-u(A)=q^-_n$. This shows $h_{n+1}= h_n$.

Let $(A,B)$ be one of the $h_n = h_{n+1}$ nonintersecting pairs of subsets of $[n+1]$ with $u(B)-u(A)=q^-_{n+1}$ as above. As before, since  $q^-_{n+1}\equiv 0 \pmod4$, neither of the sets contain~$1$. We can use these pairs to construct the same number of nonintersecting pairs  of $\preceq_{2{\bf \ell}_{n+1}}$ with utility difference $q^+_{n+1}=q^-_{n+1}+2$. Indeed, adding~$1$ to $B$ will create a pair $(A,B\cup \{1\})$ with a utility difference $q^-_{n+1}+2=q^+_{n+1}$. We can also use (\ref{-4-}) and a shifting technique to create another $g_n$ nonintersecting pairs with utility difference $q^+_{n+1}$. Indeed, if $(A,B)$ is one of the $g_n$ nonintersecting pairs in $\preceq_{2{\bf \ell}_{n}}$ with utility difference $q^+_{n}$, then the pair $(\{1\}\cup \overline{A}, \overline{B})$ will be nonintersecting in $\preceq_{2{\bf \ell}_{n+1}}$ with utility difference $2q^+_n-2=q^+_{n+1}$. We observe that the $h_{n+1}$ pairs $(C,D)$ constructed in the first method all have $1\in D$ while the $g_n$ pairs $(C,D)$ constructed in the second method all have $1 \in C$, and so the two methods never construct the same pair. Thus $g_{n+1}\ge g_n+h_n$. 

Now, let $(C,D)$ be any nonintersecting pair in $\preceq_{2{\bf \ell}_{n+1}}$ with  utility difference $u(D)-u(C)=q^+_{n+1}$. As $n+1$ is even,  Proposition~\ref{pr6} gives $q^+_{n+1}\equiv 2\pmod4$. This implies that either $1\in C$ or $1\in D$. Now as above, we can show that $(C,D)$ can be obtained as $(\{1\}\cup \overline{A},\overline{B})$ or $(A,\{1\}\cup B)$ by the second or the first method, respectively. Thus $g_{n+1}= g_n+h_n$. 

For even $n$, the statement can be proved similarly, using the other two equations in Proposition~\ref{pr7} and congruences in Proposition~\ref{pr6}.
\end{proof}

\noindent{\it Proof of Theorem~\ref{Dtheorem} (a).}
Let us consider the case $n=3$. We have $q_3=3$, so $q_3^-=2$ and $q_3^+=4$. We have three nonintersecting pairs in $\preceq_{2{\bf \ell}_3}$ with utility difference two, namely $\emptyset \preceq_{2{\bf \ell}_3} \{1\}$, $\{1\}\preceq_{2{\bf \ell}_3} \{2\}$, and $\{1,2\}\preceq_{2{\bf \ell}_3} \{3\}$, and two nonintersecting pairs with utility difference four, namely, $\emptyset \preceq_{2{\bf \ell}_3} \{2\}$ and $\{2\}\preceq_{2{\bf \ell}_3} \{3\}$. Thus $g_3=2$ and $h_3=3$. Alternatively, we may say that $(g_3, h_3)=(F_3, F_4)$. It is also easy to check that $(g_4, h_4)=(5,3)=(F_5, F_4)$.  A simple induction argument with the use of Lemma~\ref{ghrecrel} now shows that $(g_n,h_n)=(F_n,F_{n+1})$ for odd $n$ and $(g_n,h_n)=(F_{n+1},F_n)$ for even $n$.
By Lemma~\ref{g_n+h_n} we find that the number of flippable pairs of $\preceq_{\bf m}$ is
\[
g_n+h_n=F_{n+1}+F_n=F_{n+2}.
\]
It remains to notice that $\preceq_{\bf m}$ is in ${\mathcal G}_{n+1}$.
\par\medskip


\noindent{\it Proof of Theorem~\ref{Dtheorem} (b).}
Let $\preceq$ be obtained from $\preceq_{\bf m}$ by a flip. Assume it was the pair $B \prec_{\bf m} A$ in $\preceq_{\bf m}$ which was flipped, so in $\preceq$ we have $A \prec B$. Assume $\preceq$ is not representable. By Theorem~\ref{KPStheorem} there must exist a trading transform  $(\row Ak; \row Bk)$ such that $A_i\prec B_i$ for $i=1,\ldots, k$. For each $i$ we may assume that $A_i \cap B_i = \emptyset$ since otherwise we could remove the intersection for each pair and obtain another trading transform with empty intersections.

Since each element of $[n]$ appears in the sequence $\row Ak$ exactly as many times as in $\row Bk$, for the weight function $u$ of $\preceq_{\bf m}$, we must have $\sum_{i=1}^k u(A_i) = \sum_{i=1}^k u(B_i)$. However the only nonintersecting pair $C \prec D$ in $\preceq$ with $u(C) \geq u(D)$ is the flipped pair $A \prec B$, and furthermore we know from Lemma~\ref{lemma_abc} that $u(A) - u(B) =1$ and for every other pair $(A_i,B_i)$ different from $(A,B)$, $u(A_i) - u(B_i) \leq -1$. Hence for  $\sum_{i=1}^k u(A_i) = \sum_{i=1}^k u(B_i)$ to hold at least half of the pairs $A_i \prec B_i$ must be the pair $A \prec B$. Without loss of generality assume that $A_i \prec B_i$ is the pair $A \prec B$ for $i = 1,2,...,r$ with $r \geq \frac{k}{2}$.

Let $j \in A$ be any element of $A$. Then $j$ appears $r$ times in the the sequence $\row Ar$ and no times in the sequence $\row Br$. Hence it must appear $r$ times in $(B_{r+1},\ldots, B_k)$, but $r \geq \frac{k}{2}$ and $j$ can appear at most once in each $B_i$ and so we must have $r=\frac{k}{2}$, $j \in B_i$ and $u(A_i)-u(B_i) = -1$ for $i=r+1,...,k$. But $j$ was an arbitrary element of $A$, so $A \subseteq B_i$ for $i=r+1,...,k$. The same argument shows that $B \subseteq A_i$ for $i=r+1,...,k$. But if $A_i = B \cup C_i$ and $B_i = A \cup D_i$ with $C_i$, $D_i$, $A$ and $B$ all disjoint for $i=r+1,...,k$ then
\[
u(A_i)-u(B_i) =u(B \cup C_i) - u(A \cup D_i) = u(B) + u(C_i) - u(A) - u(D_i) = -1.
\]
But $u(A) - u(B) = 1$ and so $u(C_i) = u(D_i)$. Since $\preceq_{\bf m}$ is a linear order, this implies $C_i = D_i = \emptyset$ and so $B \prec  A$ which gives the desired contradiction.


\section{The Upper Bound}

We now present a result giving an upper bound on the number of flippable pairs in any comparative probability order, representable or not. This will give us an upper bound for $m(n)$ and hence for $M(n)$. The basic result is in the following lemma which estimates the number of flippable pairs from above.

\begin{lemma}
\label{ub}
Let $\preceq$ be a comparative probability order on $n$ atoms. If $s$ is any positive integer such that $\sum_{i=0}^s2^{i}\binom{n}{i}\geq2^{n}-1$,
then $|\text{Irr}({\cal C(\preceq)})| \le \sum_{i=0}^s\binom{n}{i}$.
\end{lemma}

\begin{proof}
We first prove that if $A\prec B$ and $E\prec F$ are two distinct flippable pairs then $A\cup B\neq E\cup F$. Let $A\prec B$ be a flippable pair and consider $\preceq$ restricted to the subsets of $D=A\cup B$ and call this order $\preceq^{\prime}$. Clearly $\preceq^{\prime}$ is a comparative probability order:
\[
\emptyset=D_{1}\prec^{\prime}D_{2}\prec^{\prime}D_{3}\prec^{\prime}\ldots\prec^{\prime}D_{2^{r}-1}\prec^{\prime}D_{2^{r}}=D
\]
where $r=|D|$, $D_{i}\subseteq D$ and $D_{i}\prec^{\prime}D_{j}\Longleftrightarrow D_{i}\prec D_{j}$.
Because $A$ and $B$ were adjacent in $\preceq$, they will also be adjacent in $\preceq^{\prime}$, and since $A$ and $B$ are complements in $D$, they must be the central pair of $\preceq^{\prime}$, i.e. $A=D_{2^{r-1}}$ and $B=D_{2^{r-1}+1}$. However if $E\prec F$ was also a flippable pair with $E\cup F=D$, then it must also be the central pair of $\preceq^{\prime}$, and hence $A=E$ and $B=F$. 

We now look at $\preceq$:
\[
\emptyset=A_{1}\prec A_{2}\prec A_{3}\prec\ldots\prec A_{2^{n}-1}\prec A_{2^{n}}=[n].
\]

Call the gap between two adjacent subsets an {\em adjacency}. There are a total of $2^{n}-1$ adjacencies, one for each $\prec$ sign in the order above. Consider a flippable pair $A\prec B$ in $\preceq$ and let $r=|(A\cup B)^c|$, which the size of the complement of $A\cup B$. From the definition of flippable pairs, every pair of the form $A\cup C\prec B\cup C$, where $C\subseteq (A\cup B)^c$, is adjacent. Let $C_{1},C_{2},\ldots,C_{2^{r}}$ be the subsets of $(A\cup B)^c$. Every pair $A\cup C_{i}\prec B\cup C_{i}$ will take up an adjacency, and hence the flippable pair $A\prec B$ will take up exactly $2^{r}$ adjacencies.

Hence we know that for every $r$ at most $\binom{n}{r}$ flippable pairs take up exactly $2^{r}$ adjacencies. This is because if there were more than $\binom{n}{r}$ such flippable pairs then by the pigeonhole principle two of the flippable pairs $A_{1}\prec B_{1}$ and $A_{2}\prec B_{2}$ will have $(A_{1}\cup B_{1})^c=(A_{2}\cup B_{2})^c$ and so must be the same pair. Hence there can be at most $\binom{n}{0}=1$ flippable pair that takes up $1$ adjacency, at most $\binom{n}{1}$ flippable pairs that take up $2^1$ adjacencies, at most $\binom{n}{2}$ flippable pairs that take up $2^2$ adjacencies, etc. But we have only $2^{n}-1$ adjacencies, so if we choose $s$ such that $\sum_{i=0}^s2^{i}\binom{n}{i}\geq2^{n}-1$ then  the number of flippable pairs cannot exceed $\sum_{i=0}^s \binom{n}{i}$.
\end{proof}

While the result is true for any such $s$, to maximize the strength of the upper bound we clearly wish to take the smallest value of $s$ possible. We will further need the binary entropy function ${H(\lambda)=-\lambda\log\lambda-(1-\lambda)\log(1-\lambda)}$ where the logarithms are of base 2. Using known approximations to the binomial coefficient we obtain the following:

\begin{corollary}
\label{c1}
Let $\lambda$ be the solution to the equation $\lambda+H(\lambda)=1$ and $\lambda<c<\frac{1}{2}$. Then 
$m(n)\leq2^{H(c)n}$ for sufficiently large $n$. In particular, for any $\lambda<c<\frac{1}{2}$ we have $m(n)=O\left(2^{H(c)n}\right)$.
\end{corollary}

\begin{proof}
Let $\lambda<c<\frac{1}{2}$ and $\lambda < c' < c$ (e.g. $c' = \frac{\lambda + c}{2}$). Then $c'+H(c')>1$ since $H(x)$ is strictly increasing for $0<x<\frac{1}{2}$, and for sufficiently large $n$ it holds that $\lfloor cn \rfloor > \lceil c'n \rceil$. Hence
\[
\underset{i=0}{\overset{\lfloor cn \rfloor}{\sum}}2^{i}\binom{n}{i} > 2^{\lceil c'n \rceil}\binom{n}{\lceil c'n \rceil} \geq 2^{c'n}\frac{\sqrt{\pi}}{2}\frac{1}{\sqrt{2\pi nc'(1-c')}} 2^{H(c')n} > 2^n
\]
where the second inequality is obtained from \cite[p. 466]{PW} and the last inequality holds for sufficiently large $n$. So by Lemma~\ref{ub}, we have for sufficiently large $n$
\[
m(n)\leq\underset{i=0}{\overset{\lfloor cn \rfloor}{\sum}}\binom{n}{i}\leq c^{-cn}(1-c)^{-(1-c)n}=2^{H(c)n}
\]
where the second inequality is also obtained from \cite[p. 468]{PW}.
\end{proof}

\begin{example}
\label{upeg}
Take $c=0.25$ and consider $s=\lfloor{cn}\rfloor$. It can be checked that for $n \geq 102$ and $c' = c - \frac{1}{102} \leq \frac{s}{n}$ we have $\lceil c'n \rceil \leq s$ and the following inequalities hold
\[
\underset{i=0}{\overset{s}{\sum}}2^{i}\binom{n}{i} > 2^{\lceil c'n \rceil}\binom{n}{\lceil c'n \rceil} \geq 2^{c'n}\frac{\sqrt{\pi}}{2}\frac{1}{\sqrt{2\pi n c' (1-c')}}2^{H(c')n} > 2^n.
\]
Hence by Lemma~\ref{ub} it holds that
\[
m(n)\leq\underset{i=0}{\overset{\lfloor cn \rfloor}{\sum}}\binom{n}{i}\leq 2^{H(c)n}.
\]
Here $2^{H(c)} < 1.7548$.
Along with Theorem~\ref{Dtheorem} and standard bounds on the Fibonacci sequence, we have the following bounds for $n \geq 102$:
\[
F_{n+1} = \left[ \frac{\phi^{n+1}}{\sqrt{5}}\right] \leq m(n) \leq 1.7548^n,
\]
where $\phi \approx 1.6180$ is the golden ratio and $[x]$ is the closest integer to $x$.
\end{example}

Clearly in this example the exponent of 2 in the upper bound of $m(n)$ can be brought arbitrarily close to $H(\lambda)$ for sufficiently large $n$. As $2^{H(\lambda)} \approx 1.7087$, this gives the rough bounds $1.6180^n < M(n) \leq m(n) < 1.7087^n$ up to constant factors.

\begin{corollary}
For sufficiently large $n$
\[
1.6180^n < M(n) < 1.7087^n
\]
up to constant factors.
\end{corollary}


\section{Further Research}

We would like to know, of course, if Searles' conjecture is  true. Or at least, we would like to reduce the gap between the current bounds for $M(n)$ further.
There are some interesting questions that are not directly related to  Searles' conjecture but nevertheless  interesting.  One of them is the question of connectedness of 
${\cal G}_n$. Since the subgraph of representable comparative probability orders is clearly connected, this question shows that we do not really understand much about 
nonrepresentable orders. In particular, this question will be answered in the affirmative if we could show that any order is connected to a representable order by a series of flips. A similar 
question is to find the minimum value of  $ |{\rm Irr}({\cal C})|$ in ${\cal G}_n$. For representable orders this minimum is $n$ but for nonrepresentable orders we cannot even say if it is possible that the central pair is the only flippable pair in the order. Maclagan also emphasised these questions \cite{DM}.


\end{document}